# High order elements in extensions of finite fields given by binomials

Roman Popovych
Department of Specialized Computer Systems, Lviv Polytechnic National University,
Lviv, Ukraine
e-mail: rombp07@gmail.com

**Abstract**

The extension of finite field $F_q$, given by a binomial, is of the form $F_q(\theta) = F_{q^m} = F_q[x]/(x^m - a)$. We consider any such extension and construct in it elements with high multiplicative order. A general method is described, which is an improvement of our method (Finite Fields Appl., 19 (1), 86–92, 2013). We take element $\theta + b$ ($b \in F_q^*$), which is a linear binomial in variable $\theta$. Consecutively raising the binomial to a power, we obtain more linear binomials. Then we form certain non-linear binomials from each linear binomial. In total, the initial element generates $m$ binomials. We construct pairwise distinct products of these binomials and obtain a lower bound on the number of the products, which is a lower bound $2^{\sqrt{2m}}$ on the order of $\theta + b$. This result improves the best previously known bound $5^{\sqrt[3]{m/2}}$. A numerical example that illustrates the approach is provided as well.

It is well known that the multiplicative group of a finite field is cyclic. A generator of the group is called primitive element. The problem of constructing efficiently a primitive element for a given finite field is notoriously difficult in the computational theory of finite fields. That is why one considers less restrictive question: to find an element with high multiplicative order. We are not required to compute the exact order of the element. It is sufficient in this case to obtain a lower bound on the order. High order elements are needed in several applications. Such applications include but are not limited to cryptography, coding theory, pseudo random number generation and combinatorics.

**Previous work.** The extension specified by a binomial is of the form $F_q[x]/(x^m - a)$. It is shown in [4] how to construct high order element in such extension with the condition that $m$ divides $q-1$. The lower bound $5{,}8^m$ is obtained in this case. High order elements are constructed in [3, 5] for extensions $F_{q^m}$ ($m = 2^t$, $q \equiv 1 (\mathrm{mod}\, 4)$, lower bound $2^{(t^2+3t)/2 + ord_2(q-1)}$) and in [3] for extensions $F_{q^m}$ ($m = 3^t$, $q \equiv 1 (\mathrm{mod}\, 3)$, $q \neq 4$, lower bound $3^{(t^2+3t)/2 + ord_3(q-1)}$) without the mentioned before division condition. These extensions are considered as recursive towers of finite fields, but can be also specified by binomials. For arbitrary $m$ and without the division condition, the best known results are: the lower bound $2^{\sqrt[3]{2m}}$ [10] and the refined bound $5^{\sqrt[3]{m/2}}$ [2]. A fairly complete list of references to works, related to the construction of high order elements in finite fields, is in [5].

**Our results.** Throughout this paper $q$, $m$ and $a$ are such that the extension $F_q(\theta) = F_{q^m} = F_q[x]/(x^m - a)$, where $\theta$ is the coset of $x$, exists. It is clear that $\theta^m = a$. We consider any extension of this form and construct in it elements with the multiplicative order at least $2^{\sqrt{2m}}$.

In [2, 10], the fact that $m = kl$ was used and two elements were constructed. The order of the first element depends on $k$, and the order of the second one depends on $l$. The idea was as follows: if $q-1$ has a big divisor $k$, we use for the construction the method from [4]; if $q-1$ has no a big divisor $k$, then $l$ is big, and we use for the construction the method similar to that in [1, 9, 11, 12]. As the product $m = kl$ is fixed, then such approach gave mentioned before results from [2, 10].

In this paper an improvement of the method from [2] is suggested. We take only one element, which can be considered as linear binomial in variable $\theta$. Consecutively raising the binomial to the power $q^l$, we obtain $k-1$ more linear binomials. From each linear binomial we form $l-1$ non-linear binomials of degrees $ik+1\,(i = 1,2,...,l-1)$. Then we construct pairwise distinct products of all $kl$ binomials. A lower bound on the number of these products gives the result $2^{\sqrt{2m}}$.

Our main result is the following theorem.

**Theorem 1.** *Let $b$ be any non-zero element in $F_q$. Then element $\theta + b$ of the field $F_q(\theta) = F_q[x]/(x^m - a)$ has the multiplicative order at least $2^{\sqrt{2m}}$.*

This result improves two best known results: the lower bound $2^{\sqrt[3]{2m}} = 2{,}39^{\sqrt[3]{m}}$ [10] and the refined bound $5^{\sqrt[3]{m/2}} = 3{,}58^{\sqrt[3]{m}}$ [2].

## 1. Preliminaries

Throughout this paper $F_q$ is a field of $q$ elements, where $q$ is a power of a prime number, and $F_q^*$ denotes the multiplicative group of the field. For positive integers $v$ and $w$, $\mathrm{ord}_w v$ means the order of $v$ modulo $w$. For an integer $n$, $Z_n^*$ is the multiplicative group of integers modulo $n$. By $|G|$ we denote the number of elements in the set $G$.

The relationship between $q$ and $m$ and the choice of $a$ are well known [6–8]. In finite fields of characteristic two there is only one irreducible binomial $x-1$. In odd characteristic, we can test $x^m - a$ for irreducibility using [6, Theorem 3.75].

**Theorem 2.** *Let* $m \geq 2$ *be an integer and* $a \in F_q^*$. *Then the binomial* $x^m - a$ *is irreducible in* $F_q[x]$ *if and only if the following two conditions are satisfied:*

1. *Each prime factor of* $m$ *divides the order* $e$ *of* $a \in F_q^*$, *but not* $(q-1)/e$;
2. $q \equiv 1 (\bmod 4)$ *if* $m \equiv 0 (\bmod 4)$.

The first condition of Theorem 2 means that $\gcd((q-1)/e, m) = 1$. As a development of Theorem 2, given a number $q$, it is described precisely in [8] for which degrees $m$ there exist irreducible binomials, and also element $a$ is constructed explicitly. In more detail, two conditions are as follows: each prime factor of $m$ divides $q-1$ and if 4 divides *m*, then 4 divides $q-1$. In the case $q = 3$ the only possible extension is for $m = 2$. If $q \geq 5$ is odd, then they can construct the extensions for infinitely many $m$. Therefore, we take to the end of the paper that field characteristic is odd and $q \geq 5$.

Let $q - 1 = p_1^{e_1} \ldots p_r^{e_r}$ be the factorization in pairwise distinct primes $p_i (1 \leq i \leq r)$, then $m = p_{s_1}^{l_1} \ldots p_{s_t}^{l_t}$, where $t \leq r$ and $\{p_{s_1}, \ldots, p_{s_t}\} \subseteq \{p_1, \ldots, p_r\}$ [8]. To simplify notation, enumerate these primes in such a way that $\{p_{s_1}, \ldots, p_{s_t}\} = \{p_1, \ldots, p_t\}$ and $m = p_1^{l_1} \ldots p_t^{l_t}$. Note that if $p_i = 2$ for some $1 \leq i \leq t$ and $l_i \geq 2$, then $e_i \geq 2$. Element $a$ with the order $e = p_1^{e_1} \ldots p_t^{e_t}$ equals $\alpha^{(q-1)/e}$, where $\alpha$ is a primitive element in $F_q^*$.

The following lemma is a reformulation of [10, Lemma 3].

**Lemma 3.** *Let* $q-1$, *m and e be as described above and* $1 \leq i \leq t$. *If either* $p_i$ *is odd and* $l_i \geq 1$ *or* $p_i = 2$, $l_i \geq 2$ *and* $e_i \geq 2$, *then the order of* $q$ *modulo* $p_i^{l_i}$ *equals* $\tau(p_i^{l_i - e_i}) = \begin{cases} 1, \text{if } l_i \leq e_i \\ p_i^{l_i - e_i}, \text{if } l_i > e_i \end{cases}$. *If* $p_i = 2$ *and* $l_i = 1$, *then the order of* $q$ *modulo* $p_i^{l_i}$ *equals* 1.

Applying Lemma 3 and the Chinese remainder theorem, we obtain that the order $l = ord_m q$ is equal to

$$l = \prod_{1 \leq i \leq t} ord_{p_i^{l_i}} q = \prod_{1 \leq i \leq t} \tau(p_i^{l_i - e_i}) = \prod_{\substack{1 \leq i \leq t \\ l_i > e_i}} p_i^{l_i - e_i} . \tag{1}$$

Relation (1) implies that if $m' = p_i m$ for some divisor $p_i$ of *l*, then $ord_{m'} q = p_i ord_m q$. The obvious generalization of this fact is as follows.

**Lemma 4.** *Let* $q-1$, $m$ *and* $l$ *be as described above.* Assume that $m' = mk'$, where $k' > 1$ is a divisor of *l*. Then $ord_{m'} q = k' \cdot ord_m q$.

We also have $m = kl$, where $k = m/l = \prod_{1\le i\le t} p_i^{l_i} / \tau(p_i^{l_i - e_i})$. Clearly if $l_i \le e_i$, then $p_i^{l_i} / \tau(p_i^{l_i - e_i}) = p_i^{l_i}$, and if $l_i > e_i$, then $p_i^{l_i} / \tau(p_i^{l_i - e_i}) = p_i^{e_i}$. Therefore $k = \prod_{1\le i\le t} p_i^{\min(l_i, e_i)} = \gcd(q-1, m)$.

Hence, $k$ is a divisor of $e$ (obviously $e$ divides $q-1$), $k = \gcd(q-1, m)$ and $\gcd((q-1)/k, l) = 1$. Most of the mentioned facts are summarized in the following lemma which is a slight modification of [10, lemma 4].

**Lemma 5.** *Let $l$ be the order of $q$ modulo $m$. Then $m = kl$, where $k$* is a divisor of *e*, $k = \gcd(q-1, m)$, $l$ is coprime with $(q-1)/k$ *and the subgroup $\langle q \rangle$ of $Z_m^*$ can be written in the form* $\langle q \rangle = \{ik + 1 \mid i = 0, ..., l-1\}$.

Obviously, the condition that $m$ divides $q-1$ is equivalent to $l = 1$ ($m = k$). Construction of high order elements in this case was considered in [4]. So, we assume in this paper that $m$ does not divide $q-1$. Then possible values are $l \ge 2$ and $k \ge 3$.

As *l* is the order of *q* modulo *m*, then $q^l = 1 \bmod m$, that is

$$q^l = 1 + Tm \tag{2}$$

for some integer *T*. Clearly $a^T = a^{T \bmod (q-1)}$. Given below Lemma 6 describes a feature of the integer.

**Lemma 6.** *Let k be as in Lemma 5. Then the order of element $a^T$ in the multiplicative group $F_q^*$ equals k.*

**Proof.** Using the well known fact about the order of element $a^T$, we have:

$$ord(a^T) = \frac{ord(a)}{\gcd(ord(a), T)} = \frac{e}{\gcd(e, T)}. \tag{3}$$

Recall that *k* divides *e*, i. e. $e = fk$ for some integer *f*, and that $m = kl$, where $k$ is a divisor of $q-1$ (see Lemma 5). Equality (2) implies $Tl = \frac{q-1}{k}(q^{l-1} + ... + q + 1)$. As, according to Lemma 5, $l$ is coprime with $\frac{q-1}{k}$, then $q^{l-1} + ... + q + 1$ is divisible by $l$ and $T = \frac{q-1}{k} \cdot \frac{q^{l-1} + ... + q + 1}{l}$. Then we can write

$$\gcd(e, T) = \gcd(fk, f\frac{q-1}{e} \cdot \frac{q^{l-1} + ... + q + 1}{l}) = f \gcd(k, \frac{q-1}{e} \cdot \frac{q^{l-1} + ... + q + 1}{l}).$$

Since, according to Theorem 2, $\gcd(m, \frac{q-1}{e}) = 1$, then $\gcd(k, \frac{q-1}{e}) = 1$ and

$$\gcd(e,T) = f\,\gcd(k, \frac{q^{l-1}+...+q+1}{l}).$$

Denote $k' = \gcd(k, \frac{q^{l-1}+...+q+1}{l})$, $m' = mk'$ and assume that $k' > 1$. As $k$ divides $q$-1 and $\frac{q^{l-1}+...+q+1}{l}$ divides $q^{l-1}+...+q+1$, then $k'$ is a common divisor of $q-1$ and $q^{l-1}+...+q+1$. Therefore $q \equiv 1 \bmod k'$ and $q^{l-1}+...+q+1 \equiv l \bmod k'$. Hence $k'$ divides $l$.

Clearly $lk'$ divides $q^{l-1}+...+q+1$, $m' = klk'$ divides $q^l - 1$ and so $q^l \equiv 1 \bmod m'$. We claim that $l = ord_{m'}q$. Indeed, if there exists an integer $l' < l$ such that $q^{l'} \equiv 1 \bmod m'$, then $q^{l'} \equiv 1 \bmod m$ - a contradiction to the fact $l = ord_m q$. Hence we have $l = ord_{m'}q = ord_m q$ - a contradiction to Lemma 4. So, $k' = 1$ and $\gcd(e,T) = f$. Taking into account equality (3), we obtain $ord(a^T) = \frac{e}{f} = k$. □

## 2. Explicit construction of high order elements

Below we construct explicitly in the field $F_q(\theta) = F_q[x]/(x^m - a)$ elements with the order at least $2^{\sqrt{2m}}$.

**Theorem 7.** *Let $m,k,l$ be as in Lemma 5 and $b$ be any non-zero element in $F_q$. The subgroup, generated by element $\theta + b$, contains the following pairwise distinct elements:*

$$a^{jT+r_i}\theta^{ik+1} + b, \tag{4}$$

where $0 \le i \le l-1$, $0 \le j \le k-1$ and $r_i$ are some integers (in particular, $r_0 = 0$).

**Proof.** Using equality (2), we can write $\theta^{q^l} = \theta^{Tm+1} = (\theta^m)^T\theta = a^T\theta$. Based on this fact, we successively raise element $\theta + b$ (linear binomial in $\theta$) to the power $q^l$:

$$(\theta+b)^{q^l} = a^T\theta + b,\ (a\theta+b)^{q^l} = a^{2T}\theta + b, \ldots,\ (a^{(k-2)T}\theta+b)^{q^l} = a^{(k-1)T}\theta + b.$$

As the order of $a^T$ equals $k$ (see Lemma 6), we obtain $k$ pairwise distinct linear binomials $a^{jT}\theta + b$ ($j = 0,1,...,k-1$).

According to Lemma 5, for each $i = 0,1,...,l-1$, an integer $\alpha_i \in \{0,...,l-1\}$ exists such that $q^{\alpha_i} \equiv (ik+1) \bmod m$, that is $q^{\alpha_i} \equiv (ik+1) + r_i \cdot m$ for some integer $r_i$. Then $\theta^{q^{\alpha_i}} = \theta^{i \cdot k+1}(\theta^m)^{r_i} = a^{r_i}\theta^{ik+1}$. Obviously $\alpha_0 = 0$ and $r_0 = 0$. So raising every linear binomial $a^{jT}\theta + b$ ($j = 0,1,...,k-1$) to powers $q^{\alpha_i}$ ($i = 1,...,l-1$), we obtain $l-1$ non-linear binomials:

$$(a^{jT}\theta+b)^{q^{\alpha_i}} = a^{jT+r_i}\theta^{ik+1}+b.$$

Clearly all $kl$ binomials of the form (4) belong to the subgroup generated by element $\theta+b$. If $0\le j_1<j_2\le k-1$, then binomials $a^{j_1T+r_i}\theta^{ik+1}+b$ and $a^{j_2T+r_i}\theta^{ik+1}+b$ of the same degree $i$ are pairwise distinct. Indeed, if $a^{j_1T+r_i}=a^{j_2T+r_i}$, then $(a^T)^{(j_1-j_2)}=1$, where $j_2-j_1<k$ – a contradiction to Lemma 6. Hence, all $kl$ binomials are pairwise distinct. □

We give below an example to Theorem 7.

**Example.** Let $q$=7, then $q-1=2\cdot 3$ and we can take $m=2\cdot 3^3=54$. Take primitive element $\alpha=3$ of $F_7$ and construct $a=\alpha^{(q-1)/e}=\alpha^{6/6}=\alpha$ with the order 6. So, we have the extension $F_7(\theta)=F_7[x]/(x^{54}-a)$ and $\theta^{54}=a$.

Powers of $q$=7 modulo $m$=54 are as follows: $7^0$=1, $7^1$=7, $7^2$=49, $7^3$=19, $7^4$=25, $7^5$=13, $7^6$=37, $7^7$=43, $7^8$=31. We have $7^9=40353607=1+747289\cdot 54\equiv 1 \bmod 54$. Then $l$=9 is the order of $q$=7 modulo $m$=54 and $k=6$, $T=747289$ and $T \bmod (q-1)=1$. We also have $a^T=a^{T \bmod (q-1)}=a$ and $ord(a^{T \bmod (q-1)})=ord(a)=6=k$.

Consider the binomial $\theta+1$. Linear binomials, obtained from this binomial, are as follows:

$(\theta+1)^{7^9}=a\theta+1$, $(a\theta+1)^{7^9}=a^2\theta+1$, $(a^2\theta+1)^{7^9}=a^3\theta+1$, $(a^3\theta+1)^{7^9}=a^4\theta+1$, $(a^4\theta+1)^{7^9}=a^5\theta+1$.

Non-linear binomials, formed from $\theta+1$, are as follows:

$(\theta+1)^7=\theta^7+1$,

$(\theta^7+1)^7=\theta^{49}+1$,

$(\theta^{49}+1)^7=\theta^{343}+1=(\theta^{54})^6\theta^{19}+1=a^6\theta^{19}+1=\theta^{19}+1$,

$(\theta^{19}+1)^7=\theta^{133}+1=(\theta^{54})^2\theta^{25}+1=a^2\theta^{25}+1$,

$(a^2\theta^{25}+1)^7=a^{14}\theta^{175}+1=a^2(\theta^{54})^3\theta^{13}+1=a^5\theta^{13}+1$,

$(a^5\theta^{13}+1)^7=a^{35}\theta^{91}+1=a^5\theta^{54}\theta^{37}+1=a^6\theta^{37}+1=\theta^{37}+1$,

$(\theta^{37}+1)^7=\theta^{259}+1=(\theta^{54})^4\theta^{43}+1=a^4\theta^{43}+1$,

$(a^4\theta^{43}+1)^7=a^{28}\theta^{301}+1=a^4(\theta^{54})^5\theta^{31}+1=a^3\theta^{31}+1$.

Below we write other obtained non-linear binomials, but do not give the correspondent calculations, because they are analogous to the given above calculations:

$a\theta+1$, $a\theta^7+1$, $a\theta^{49}+1$, $a\theta^{19}+1$, $a^3\theta^{25}+1$, $\theta^{13}+1$, $a\theta^{37}+1$, $a^5\theta^{43}+1$, $a^4\theta^{31}+1$;

$a^2\theta+1$, $a^2\theta^7+1$, $a^2\theta^{49}+1$, $a^2\theta^{19}+1$, $a^4\theta^{25}+1$, $a\theta^{13}+1$, $a^2\theta^{37}+1$, $\theta^{43}+1$, $a^5\theta^{31}+1$;

$a^3\theta+1$, $a^3\theta^7+1$, $a^3\theta^{49}+1$, $a^3\theta^{19}+1$, $a^5\theta^{25}+1$, $a^2\theta^{13}+1$, $a^3\theta^{37}+1$, $a\theta^{43}+1$, $\theta^{31}+1$;

$a^4\theta+1$, $a^4\theta^7+1$, $a^4\theta^{49}+1$, $a^4\theta^{19}+1$, $\theta^{25}+1$, $a^3\theta^{13}+1$, $a^4\theta^{37}+1$, $a^2\theta^{43}+1$, $a\theta^{31}+1$;
$a^5\theta+1$, $a^5\theta^7+1$, $a^5\theta^{49}+1$, $a^5\theta^{19}+1$, $a\theta^{25}+1$, $a^4\theta^{13}+1$, $a^5\theta^{37}+1$, $a^3\theta^{43}+1$, $a^2\theta^{31}+1$.
In total, we have $m = kl = 54$ binomials.

Take $k$ linear binomials $a^{jT}\theta + b$ ($0 \le j \le k-1$), that, according to Theorem 7, are in the subgroup generated by element $\theta + b$. Obviously products of these binomials, in which we take every binomial at most once, are pairwise distinct. If this is not the case, then we have that $\theta$ is a root of non-zero polynomial of degree smaller than *m*, which gives a contradiction Therefore we obtain the lower bound $2^k$ on the order of $\theta + b$. Considering products of both positive and negative powers of these binomials, we obtain a bit better lower bound $5{,}8^k$. The technique is well known and described, in particular, in [4]. Hence the following lemma is true.

**Lemma 8.** *Let* $m, k$ *be as in Lemma 4 and* $b$ *be any non-zero element in* $F_q$. *Then* $\theta + b$ *has in the field* $F_q(\theta)= F_q[x]/(x^m - a)$ *the multiplicative order at least* $5{,}8^k$.

**Lemma 9.** *Let* $k, l$ *be as in Lemma 5. Set*. For an integer *c* ($0 \le c < l-1$), if $u_c, u_{c+1}, \ldots, u_{l-1}$, $v_c$ *are non-negative integers and* $u_c \le k$, *then the equality*

$$u_c(ck+1) = v_c(ck+1) + u_{c+1}[(c+1)k+1]\ldots + u_{l-1}[(l-1)k+1], \tag{5}$$

*holds only when* $u_c = v_c, u_{c+1} = u_{c+2} = \ldots = u_{l-1} = 0$.

**Proof.** Since $ck+1 < (c+1)k+1 < \ldots < (l-1)k+1$, then $u_c > v_c + u_{c+1} + \ldots + u_{l-1}$ and $u_c \ge v_c$. Suppose that $u_c > v_c$. Note that every number $ik+1 (c \le i \le l-1)$ equals 1 modulo $k$. So, the left side of equality (5) equals $u_c$ modulo *k*, and the right side equals $v_c + u_{c+1} + \ldots + u_{l-1}$ modulo *k*.

If $u_c < k$, then numbers $u_c$ and $v_c + u_{c+1} + \ldots + u_{l-1}$ are different modulo *k*. Hence, equality (5) does not hold. If $u_c = k$, then $u_c = 0 \bmod k$. However, $v_c + u_{c+1} + \ldots + u_{l-1}$ is less than $u_c$ and does not equal 0 modulo *k*. Thus equality (5) does not hold as well.

Hence, the only possible case is $u_c = v_c, u_{c+1} = u_{c+2} = \ldots = u_{l-1} = 0$. □

We construct below products of degree smaller than *m* of described in Theorem 7 binomials, taking every binomial at most once. Let us consider the set *S* of solutions $(e_{ij})_{0\le i\le l-1, 0\le j\le k-1}$ of the linear Diophantine inequality

$$\sum_{0\le i\le l-1} (ik+1) \sum_{0\le j\le k-1} e_{ij} < m, \tag{6}$$

for which $e_{ij} \in \{0,1\}$.

**Theorem 10.** *Let* $m,k,l$ *and* $d_i (i=0,...,l-1)$ *be as in Lemma 5 and* $b$ *be any non-zero element in* $F_q$. *Then* $\theta+b$ *has in* $F_q(\theta)=F_q[x]/(x^m-a)$ *the multiplicative order at least the number of elements in the set* $S$.

**Proof.** For every element $(e_{ij})$ from *S* we construct the following product

$$\prod_{0\le i\le l-1}\prod_{0\le j\le k-1}(a^{jT+r_i}\theta^{ik+1}+b)^{e_{ij}},$$

which, according to Theorem 7, belongs to the subgroup generated by $\theta+b$. Clearly it suffices to show that if two elements from *S* are distinct, then the corresponding products are not equal.

Assume that elements $(e_{ij})$ and $(f_{ij})$ from $S$ are distinct, and the corresponding products are equal:

$$\prod_{0\le i\le l-1}\prod_{0\le j\le k-1}(a^{jT+r_i}\theta^{ik+1}+b)^{e_{ij}}=\prod_{0\le i\le l-1}\prod_{0\le j\le k-1}(a^{jT+r_i}\theta^{ik+1}+b)^{f_{ij}}.$$

Since $x^m-a$ is the characteristic polynomial of $\theta$, we write

$$\prod_{0\le i\le l-1}\prod_{0\le j\le k-1}(a^{jT+r_i}x^{ik+1}+b)^{e_{ij}}=\prod_{0\le i\le l-1}\prod_{0\le j\le k-1}(a^{jT+r_i}x^{ik+1}+b)^{f_{ij}}(\mathrm{mod}(x^m-a)).$$

As, according to the definition of the set *S*, there are polynomials of degree smaller than $m$ on the left and on the right side of the equality, these polynomials are equal as polynomials over $F_q$:

$$\prod_{0\le i\le l-1}\prod_{0\le j\le k-1}(a^{jT+r_i}x^{ik+1}+b)^{e_{ij}}=\prod_{0\le i\le l-1}\prod_{0\le j\le k-1}(a^{jT+r_i}x^{ik+1}+b)^{f_{ij}}.$$

For $i=0,...,l-1$ denote $G_i=\{j\mid j\in\{0,1...,k-1\},e_{ij}=1\}$ and $H_i=\{j\mid j\in\{0,...,k-1\},f_{ij}=1\}$. Then we can rewrite the previous equality as follows:

$$\prod_{0\le i\le l-1}\prod_{j\in G_i}(a^{jT+r_i}x^{ik+1}+b)=\prod_{0\le i\le l-1}\prod_{j\in H_i}(a^{jT+r_i}x^{ik+1}+b). \tag{7}$$

Let $c\ (0\le c\le l-1)$ be the smallest integer such that $G_c\ne H_c$. After removing common factors, related with sets $G_i=H_i\ (0\le i\le c-1)$, on both sides of (7) we obtain

$$\begin{aligned}&\prod_{j\in G_c}(a^{jT+r_c}x^{ck+1}+b)\cdot\prod_{c+1\le i\le l-1}\prod_{j\in G_i}(a^{jT+r_i}x^{ik+1}+b)=\\&=\prod_{j\in H_c}(a^{jT+r_c}x^{ck+1}+b)\cdot\prod_{c+1\le i\le l-1}\prod_{j\in H_i}(a^{jT+r_i}x^{ik+1}+b)\end{aligned} \tag{8}$$

Set $G_c^/=G_c\setminus\{G_c\cap H_c\}$ and $H_c^/=H_c\setminus\{G_c\cap H_c\}$. Clearly $G_c^/,H_c^/$ are disjoint. After removing common factors, related with the set $G_c\cap H_c$, on both sides of (8) this identity leads to the following one:

$$\prod_{j\in G_c^{/}}(a^{jT+r_c}x^{ck+1}+b)\cdot\prod_{c+1\le i\le l-1}\prod_{j\in G_i}(a^{jT+r_i}x^{ik+1}+b)=$$
$$=\prod_{j\in H_c^{/}}(a^{jT+r_c}x^{ck+1}+b)\cdot\prod_{c+1\le i\le l-1}\prod_{j\in H_i}(a^{jT+r_i}x^{ik+1}+b). \qquad (9)$$

Denote $D_1(x)=\prod_{j\in G_c^{/}}(a^{jT+r_c}x^{ck+1}+b)$ and $D_2(x)=\prod_{c+1\le i\le l-1}\prod_{j\in G_i}(a^{jT+r_i}x^{ik+1}+b)$. Note that the absolute term of polynomial $D_2(x)$ is $b^{\lambda}$, where $\lambda=\sum_{c+1\le i\le l-1}|G_i|$. We can write the left side of (9) as follows:

$$D_1(x)D_2(x)=b^{\lambda}D_1(x)+D_1(x)(D_2(x)-b^{\lambda}).$$

Polynomial $b^{\lambda}D_1(x)$ has in its expansion terms of degrees $u_c(ck+1)$, where $0\le u_c\le|G_c^{/}|\le k$, and polynomial $D_1(x)(D_2(x)-b^{\lambda})$ – terms of degrees $v_c(ck+1)+u_{c+1}[(c+1)k+1)]...+u_{l-1}[(l-1)k+1)]$, where $0\le v_c\le|G_c^{/}|\le k$, $0\le u_i\le|G_i|\le k$ for $c+1\le i\le l-1$. If $D_2(x)$ is a non-empty product, i.e. $D_2(x)\ne 1$, then at least one of numbers $u_{c+1},...,u_{l-1}$ is non-zero. Therefore according to Lemma 9, polynomial $D_1(x)(D_2(x)-b^{\lambda})$ does not have terms of degrees $u_c(ck+1)$. If $D_2(x)=1$, then clearly the left side of (9) equals $D_1(x)$. Analogous considerations can be applied to the right side of (9). Hence, we have the following equality:

$$b^{\lambda}\prod_{j\in G_c^{/}}(a^{jT+r_c}x^{ck+1}+b)=b^{\mu}\prod_{j\in H_c^{/}}(a^{jT+r_c}x^{ck+1}+b),$$

where $\mu=\sum_{c+1\le i\le l-1}|H_i|$. Introducing new variable $y=x^{ck+1}$ in the last equality, we obtain:

$$b^{\lambda-\mu}\prod_{j\in G_c^{/}}(a^{jT+r_c}y+b)=\prod_{j\in H_c^{/}}(a^{jT+r_c}y+b) \qquad (10)$$

Note that sets $G_c^{/},H_c^{/}$ are disjoint. Since $G_c\ne H_c$, then at least one of sets $G_c^{/},H_c^{/}$, is non-empty, say $G_c^{/}$. Therefore there is at least one not equal to 1 linear polynomial in variable $y$ on the left side of (10). Since $F_q[y]$ is a unique factorization ring, this polynomial must be equal to some (if any) linear polynomial on the right side (accurate to a factor from $F_q^*$). But, according to Theorem 7, such polynomials are pairwise distinct (have different coefficients near $y$, but the same absolute terms), which leads to a contradiction.

So, products corresponding to distinct elements from $S$ cannot be equal, and the result follows. □

**Lemma 11.** *Let* $m,k,l$ *be as in Lemma 5 and* $S$ *be as in Theorem 10. If* $l > 2k$*, then the number of elements in the set* $S$ *is at least* $2^{\sqrt{2m}}$.

**Proof.** Note that $m/k = l$. Rewrite the left side of inequality (6), which is in the definition of the set $S$, as the sum of terms $T_j$:

$$\sum_{0\le i\le l-1}(ik+1)\sum_{0\le j\le k-1} e_{ij} = \sum_{0\le j\le k-1}\sum_{0\le i\le l-1}(ik+1)e_{ij} = \sum_{0\le j\le k-1} T_j,$$

where $T_j = \sum_{0\le i\le l-1}(ik+1)e_{ij}$ ($j = 0,1,...,k-1$). Let us ensure for each of these terms that $T_j < m/k = l$, i. e. the following inequality holds:

$$\sum_{0\le i\le l-1}(ik+1)e_{ij} < l. \qquad (11)$$

That is, we reduce inequality (6) to $k$ inequalities of the form (11).

Below we show how to find a solution $e_{ij}$ ($0 \le i \le l-1$) of (11). Let us choose the biggest integer $w$ ($w \le l-1$) such that $\sum_{0\le i\le w}(ik+1) < l$. Recall that $k > 2$. Since

$$\sum_{0\le i\le w}(ik+1) = (wk+2)(w+1)/2 < k(w+1)^2/2,$$

we choose $w$ from the equality $k(w+1)^2/2 = l$, that is $w = \sqrt{2l/k} - 1$. Clearly if to take $e_{ij} \in \{0,1\}$ for $i = 0,...,w$ and $e_{ij} = 0$ for $i = w+1,...,l-1$ we obtain a solution of (11). The number of such solutions is $2^{w+1} = 2^{\sqrt{2l/k}}$.

Combining described above solutions of inequalities (11) for different $j$, we obtain that inequality (6) has at least $\left(2^{\sqrt{2l/k}}\right)^k = 2^{\sqrt{2lk}} = 2^{\sqrt{2m}}$ solutions. □

**Remark.** We add the condition $l > 2k$ to the statement of Lemma 11 because if $l \le 2k$, then all considerations related with obtaining of solutions of inequality (11) are vacuous.

Now we are able to prove our main result.

**Proof of Theorem 1.** Consider two possible cases.

Case 1. $k \ge l/2$. As $m = kl$, then $k \ge \sqrt{m/2}$. According to Lemma 8, the order of element $\theta + b$ is at least $5{,}8^k > 4^k = 2^{\sqrt{2m}}$.

Case 2. $l > 2k$. Applying Theorem 10 and Lemma 11, we obtain the lower bound $2^{\sqrt{2m}}$. □